\documentclass[onecolumn,final,a4paper,sort&compress]{elsarticle}
\usepackage[dvipsnames]{xcolor}
\usepackage[ruled,vlined,linesnumbered]{algorithm2e}
\RestyleAlgo{ruled}

\makeatletter
\def\ps@pprintTitle{%
 \let\@oddhead\@empty
 \let\@evenhead\@empty
 \def\@oddfoot{}%
 \let\@evenfoot\@oddfoot}
\makeatother

\usepackage[utf8]{inputenc}
\usepackage{graphicx}
\usepackage{amsmath,amssymb,amsfonts,latexsym,stmaryrd}
\usepackage{algorithmic}
\usepackage{booktabs}

\usepackage[colorlinks=true,breaklinks=true,linkcolor=lightblue,citecolor=lightblue]{hyperref}
\usepackage{stackengine}
\usepackage{mathrsfs}
\usepackage{rotating} 
\usepackage{showkeys}
\usepackage{amsthm}
\usepackage{graphicx,float,epsfig,color,fancyhdr}
\usepackage{amsopn}

\usepackage[utf8]{inputenc}
\usepackage{wrapfig}

\usepackage{amsfonts}
\usepackage{amsmath}
\numberwithin{equation}{section}
\usepackage{mathtools}

\usepackage[ruled,vlined,linesnumbered]{algorithm2e}
\let\oldnl\nl
\newcommand{\nonl}{\renewcommand{\nl}{\let\nl\oldnl}}%

\usepackage[utf8]{inputenc}
\usepackage{graphicx}
\usepackage{wrapfig}
\usepackage{comment}

\usepackage{makecell}

\usepackage[ruled,vlined,linesnumbered]{algorithm2e}
\let\oldnl\nl

\newcommand{\tenprod}{\circ}
\newcommand{\ten}[1]{\mathcal{#1}}
\newcommand{\mat}[1]{\mathbf{#1}}
\renewcommand{\vec}[1]{\mathbf{#1}}

\newcommand{\PGRAPH}[1]{\noindent\textbf{#1}}
\newcommand{\HAT}[1]{\widetilde{}}

\usepackage{lipsum}

\renewcommand{\theequation}{\thesection.\arabic{equation}}

\def\bx{\mathbf{x}}

\newcommand{\vertiii}[1]{{\left\vert\kern-0.25ex\left\vert\kern-0.25ex\left\vert #1 
    \right\vert\kern-0.25ex\right\vert\kern-0.25ex\right\vert}}

\newcommand{\TTf}{\text{TT}}
\newcommand{\FDFD}{\text{FD-FD}}
\newcommand{\SPSP}{\text{SP-SP}}

\newcommand{\Mem}{\mbox{\text{Mem}}}
\newcommand{\CmpRate}{\mbox{\text{CR}}}

\newcommand{\TNs}{\text{TNs}}

\newcommand{\EOD}{\end{document}}

\newif\ifSIAM
\SIAMfalse

\newcommand{\tentt}[1]{\ten{#1}^{TT}}
\begin{document}

\begin{frontmatter}

  \title{Tensor Network Space-Time Spectral Collocation Method for
    Time Dependent Convection-Diffusion-Reaction Equations}
  
  \author[TDIV]{Dibyendu Adak}
  \author[TDIV]{Duc P. Truong}
  \author[TDIV]{Gianmarco~Manzini}
  \author[TDIV]{Kim \O. Rasmussen}
  \author[TDIV]{Boian S. Alexandrov}
  
  \address[TDIV]{Theoretical Division, Los Alamos National Laboratory,
    Los Alamos, NM 87545, USA }

\begin{abstract}
  Emerging tensor network techniques for solutions of Partial
  Differential Equations (PDEs), known for their ability to break the
  curse of dimensionality, deliver new mathematical methods for
  ultra-fast numerical solutions of high-dimensional problems.
  Here, we introduce a Tensor Train (\TTf{}) Chebyshev spectral
  collocation method, in both space and time, for solution of the time
  dependent convection-diffusion-reaction (CDR) equation with
  inhomogeneous boundary conditions, in Cartesian geometry.
  Previous methods for numerical solution of time dependent PDEs often
  use finite difference for time, and a spectral scheme for the
  spatial dimensions, which leads to slow linear convergence.
  Spectral collocation space-time methods show exponential
  convergence, however, for realistic problems they need to solve
  large four-dimensional systems.
  We overcome this difficulty by using a \TTf{} approach as its
  complexity only grows linearly with the number of dimensions.
  We show that our \TTf{} space-time Chebyshev spectral collocation method
  converges exponentially, when the solution of the CDR is smooth, and
  demonstrate that it leads to very high compression of linear
  operators from terabytes to kilobytes in \TTf{}-format, and tens of
  thousands times speedup when compared to full grid space-time
  spectral method.
  These advantages allow us to obtain the solutions at much higher
  resolutions.
\end{abstract}

\end{frontmatter}

\raggedbottom
\section{Introduction}
\label{sec1:introduction}

Solving realistic problem often require numerical solutions of high
dimensional partial differential equations (PDEs).
The discretization of these PDEs leads to a steep rise in the
computational complexity in terms of storage and number of arithmetic
operations with each added dimension, often rendering traditional
numerical approaches impractical.
This phenomenon, known as \emph{``the curse of dimensionality'',}
cf.~\cite{bellman1966dynamic}, represents a formidable barrier to
multidimensional numerical analysis, one that persists even in the era
of exascale high-performance computing.

In recent years, tensor networks (\TNs{}) \cite{bachmayr2016tensor}
have come to the forefront as a promising strategy to counteract or
circumvent the curse of dimensionality.
\TNs{} restructure high-dimensional data into networks of
lower-dimensional tensors, enabling the division of complex data into
more manageable subsets.
Tensor factorizations such as the Canonical Polyadic Decomposition,
Hierarchical Tucker decomposition \cite{kolda2009tensor}, and Tensor
Train decompostion \cite{oseledets2011tensor} are example of tensor
networks.
Originally devised within the realm of theoretical physics, these
methods now showcase their potential for accurate and efficient
numerical solutions to high-dimensional PDEs.

Through a process we will term \emph{tensorization,} tensor networks
offer an innovative means to efficiently represent the key elements of
numerical algorithms for PDEs: grid functions that capture function
values at grid nodes, and discrete operators approximating
differential operators.
This methodology has been applied successfully to a variety of
challenging equations, from the quantum mechanical to the classical
continuum
\cite{kazeev2012low,gelss2022solving,benner2021regularization,matveev2016tensor,manzini2021nonnegative,manzini2023tensor,ye2022quantum,ye2023quantized,kormann2015semi,truong2023tensor}.

The present work considers the time-dependent
convection-diffusion-reaction (CDR) equation, a model equation of
critical importance across a broad spectrum of physical and
engineering systems.
The CDR equation enables quantitative descriptions of heat transfer,
mass transfer, fluid dynamics, and chemical interactions in complex
settings from microscale processing to atmospheric modeling.
Numerical methods for the multidimensional CDR equation also
constitute a major research area given the equation’s utility spanning
such a vast range of transport phenomena critical to climate modeling,
energy systems, biomedical systems, materials synthesis, and related
domains central to technology
innovation~\cite{stocker2011introduction,bird1960we,doi1988theory}.

Classical numerical methods, like finite differences and finite
volumes, necessitate extremely fine grids for precision, leading to
voluminous linear systems that are challenging to solve, particularly
when the problem spans four dimensions (one temporal and three
spatial).
This is why spectral collocation methods are considered, due to their
greater efficiency.

The spectral collocation method approximates the solution to
convection-diffusion-reaction (CDR) equations through the use of
high-degree polynomial interpolants, e.g., Legendre and
Chebyshev.
These interpolants are evaluated at specific collocation points chosen
within the domain.
Assuming that the groundtruth solution is regular enough, one of the
most notable advantages of this approach is the
\emph{exponential convergence}, where the interpolation error
decreases exponentially as the degree of the polynomial increases.
As a computational method for solving partial differential equations
(PDEs), spectral collocation offers several significant benefits, as
highlighted by Funaro~\cite{funaro1997spectral}: $(i)$, the method
achieves a high level of precision due to its inherent exponential
convergence as the number of degrees of freedom grows; $(ii)$, it
accurately represents solutions that exhibit complex spatial
variations by employing smooth and continuous basis functions; $(iii)$
the basis functions are globally defined over the entire domain,
making the spectral collocation method particularly suitable at
handling problems with non-local interactions and dependencies.

Yet, the standard procedure for time-dependent PDEs, which couples
spectral discretization of spatial derivatives with low-order temporal
derivatives, often results in a temporal error that overshadows the
spatial precision.
Ideally, the accuracy of time integration should align with that of
the spatial spectral approximation.
Recent advancements in this area aim to overcome the time-stepping
limitations by integrating space-time spectral collocation methods,
which have shown exponential convergence in both spatial and temporal
dimensions for sufficiently smooth heat equation solutions \cite{lui2017legendre}.

Despite these advances, the curse of dimensionality remains a
significant challenge for space-time spectral methods applied to CDR
equations, as computational complexity grows exponentially with
problem size.
In this study, we explore the discretization of the CDR equation in
both space and time, employing a \emph{space-time} operator scheme
that simultaneously addresses spatial and temporal dimensions.
This approach is particularly effective for phenomena with rapid
dynamic changes where spatial and temporal variations are closely
linked \cite{marion2013classical}.
Although the space-time method traditionally faces the drawback of
increased computational cost and memory requirements, we will
demonstrate that tensor network techniques effectively overcome
these challenges.

We utilize spectral collocation discretization within this space-time
framework, expanding the PDE solution in terms of a set of global
basis functions. For the non-periodic CDR equation solution, we employ
the Chebyshev polynomials as our basis functions.
The Chebyshev collocation method requires regular grids, that is, it works
mainly on Cartesian grids that are tensor products of 1-D domain
partitions, which is suitable to our tensor network formats.

The outline of the paper is as follows.
In Section~\ref{sec:spacetimeapproach}, we review some basic concepts
for space-time collocation methods, introduce the discretization and
its matrix formulation that leads to a linear system.
In Section~\ref{sec:Tensor-Networks}, we review the tensor notations
and definitions of the \TTf{}-format, and the \TTf{}-matrix format for linear
operators.
In Section~\ref{tensorization}, we present our \TTf{} design of the
numerical solution of the CDR equation and introduce our algorithms in
tensor train format.
In Section~\ref{results}, we present our numerical results and assess
the performance of our method.
In Section~\ref{conclusions}, we offer our final remarks and discusses
possible future work.
For completeness, part of the details of the CDR tensorization
approach are in a final appendix.



\section{Mathematical Model and Numerical Discretization}
\label{methods}

\subsection{Convection-Diffusion-Reaction Equation}
\label{eq:model:problem}

We are interested in the following CDR problem: \emph{Find
$u(t,\mathbf{x})$ such that}
\begin{align}
    \frac{\partial u}{\partial t}
    - \kappa(t,\bx) \Delta u
    + \vec{b} (t,\bx) \cdot \nabla u
    + c(t,\bx)u    &= f(t,\bx) \phantom{g(t,\bx) h(\bx)}  \hspace{-1cm}\text{in} \ (0,T] \times \Omega,\nonumber\\ 
    u              &= g(t,\bx) \phantom{f(t,\bx) h(\bx)}  \hspace{-1cm}\text{on} \ [0,T] \times \partial\Omega,\label{model:prob:bdy}\\[0.2em]
    u(0,\bx )      &= h(  \bx) \phantom{f(t,\bx) g(t,\bx)}\hspace{-1cm}\text{in} \  \Omega,\nonumber
\end{align}
where $\kappa(t,\bold{x})$, $\Vec{b}(t,\bold{x})$, and $c(t,\bx)$ are
the diffusion, convection and reaction coefficients, respectively.
Furthermore, $\Omega\subset\mathbb{R}^3$ be a three dimensional, open
parallelepipied domain with boundary $\partial \Omega$; $h(\bold{x})$
is the initial condition, (IC), and $g(t,\bold{x})$ are the boundary
conditions (BC).
Here we denote the 3D vectors, such as, $\Vec{x}=(x,y,z)$, and
matrices, in bold font.

\subsection{Chebyshev collocation method}
\label{sec:spacetimeapproach}

\if
The CDR equation contains both space and time variables and the
associated differential operators.
Here, we consider the discretization in space and time simultaneously
and thus we will work with a \emph{space-time} operator scheme that
simultaneously works in both the spatial and temporal dimensions.
Using a space-time approach simplifies the formulation and
implementation of the numerical algorithm.
The space-time approach is accurate for phenomena with rapidly
changing dynamics where the interaction between spatial and temporal
variations is significant \cite{marion2013classical}.
A traditional disadvantage of the space-time approach is the extension
of the number of dimensions, which consequently increases the
computational cost and memory requirements due to the curse of
dimensionality.
However, we will demonstrate that using tensor network machinery will
overcome this problem.
Classical numerical methods to solve the CDR equation are finite
difference approaches \cite{leveque2007finite} and spectral
collocation methods \cite{funaro1997spectral}.
\fi

We expand the approximate solution of the CDR equation on the set of
orthogonal Chebyshev polynomials $T_k(x)$ = cos(k*arcos(x)) as global
basis functions.
Then, we enforce the CDR equation at a discrete set of points within
the domain, known as \emph{collocation points}, leading to a system of
equations for the unknown expansion coefficients that we can solve.
The collocation points of the Chebyshev grid are defined as in
\cite{hussaini1989spectral}, and form the so-called
\emph{collocation grid}, cf.~\cite{funaro1997spectral}.
To this end, we need expressions for the derivatives of the
approximate solution on the Chebyshev grid.
The derivative expansion coefficients with respect to the same set of
Chebyshev polynomials are determined by multiplying the matrix
representation of the partial differential operators and the expansion
coefficients of the numerical solution.
In the rest of this section, we construct the time derivative,
${\partial}\slash{\partial t}$, the gradient, $\nabla$, and the
Laplacian, $\Delta $.
To simplify such matrix representations, we modify the Chebyshev
polynomials (following Ref.~\cite{funaro1997spectral}), and construct
a new set of $N$-th degree polynomials $l_j$, $0\leq j\leq N$, with
respect to the collocation points, $x_i$, $0\leq i\leq N$, that
satisfy:
\begin{equation}
  l_j(x_i) = 
  \left \{
  \begin{aligned}
    & 1 \qquad \text{if } \ i=j, \\
    & 0 \qquad \text{if } \ i \neq j.
  \end{aligned}
  \right.
\end{equation}
Using the new polynomials $l_j$, the discrete CDR solution becomes,  
\begin{equation}
\label{Int_legen}
    u_h(x):= \sum_{j=0}^N u(x_j) l_j(x).
\end{equation}

We construct all the spatial partial differential operators using
the matrix representation of the single, one-dimensional derivative,
$\bold{S}_x$, which reads as
\begin{equation}
  \left(\frac{\partial}{\partial x}\right)_{ij} \rightarrow 
  (\bold{S}_x)_{i j}:= \frac{d}{dx}l_j(x)\Big|_{x_i}.
  \label{time_derivative}
\end{equation}
The size of the matrix $\bold{S}_x$ is $(N+1) \times (N+1)$.
We obtain the second order derivative matrix, $\bold{S}_{xx}$, by
matrix multiplication of the first order derivatives as,
\begin{equation}
  (\bold{S}_{xx})_{ij} = \sum_{s=0}^{N} (\bold{S}_x)_{is} (\bold{S}_x)_{sj}.
\end{equation}
For completeness, we refer to \cite{funaro1997spectral} for the
derivation of this expression.

\if
As a numerical method for PDEs, the Chebyshev collocation method
demonstrates \cite{funaro1997spectral}:
$(i)$ High accuracy, based on its exponential convergence rate with the
number of degrees of freedom;
$(ii)$ Accuracy for capturing complex variations of the solution, based
on the use of smooth basis functions;
$(iii)$ A global approximation, since each basis function affects the
entire domain and is suitable for problems with non-local
dependencies.
Chebyshev collocation method requires regular grids, that is, it works
mainly on Cartesian grids that are tensor products of 1-D domain
partitions, which is suitable to our tensor network formats.
\fi

\subsubsection{Matrix form of the discrete CDR equation}

First, we introduce the \emph{space-time} matrix operators: $\frac{\partial}{\partial t}
\rightarrow$ $\bold{A}_{t}$; $ \kappa(t,\bold{x})\Delta $
$\rightarrow$ $\bold{A}_{D}$; $\vec{b}(t,\bold{x}) \cdot \nabla $
$\rightarrow$ $\bold{A}_{C}$; and $c(t,\bold{x})$ $\rightarrow$
$\bold{A}_{R}$, respectively.
Upon employing these matrices, Eq.(\ref{model:prob:bdy}) results in
the following linear system,
\begin{equation}
  (\bold{A}_{t}+ \bold{A}_{D}+\bold{A}_{C}+\bold{A}_{R}) \bold{U}= \bold{F},
  \label{mat_eq}
\end{equation}
where $\bold{U}$, and $\bold{F}$ are the vectors corresponding to the
solution, $u_h(t, \bold{x})$, and the loading term $f(t, \bold{x})$,
respectively.
The matrices $\bold{A}_{t}$, $\bold{A}_{D}$, $\bold{A}_{C}$, and
$\bold{A}_{R}$ are of size $(N+1)^4 \times (N+1)^4$ for $N^{th}$ order
spectral collocation method, while $\bold{U}$, $\bold{F}$ are column
vectors of size $(N+1)^4 \times 1$, and are designed to incorporate the
boundary conditions in Eq.(\ref{model:prob:bdy}).

\subsubsection{Time discretization using finite differences and Chebyshev grids}
Here, we display two strategies to discretize the first-order
derivative with respect to the time variable.
First, we present the well-known temporal finite difference approach
with the implicit backward Euler method \cite{li2019computational},
and second we introduce a temporal spectral collocation discretization
on a temporal Chebyshev grid.
In both strategies the space operators are discretized on a spatial
Chebyshev grid.

\PGRAPH{$\bullet$ Finite difference approach}: With $N+1$ time points, $t_0,
t_1, \ldots, t_{N}$, such that $t_0=0$ and $t_{N}=T$ are the initial
and final time points, the length of the time step is $\Delta
t=(T-T_{0})/N$.
We emphasize that the backward Euler method is unconditionally stable
and thus the stability is independent of the size of the time-step
$\Delta t$ \cite{li2019computational}.
To consider the finite differences approach, we need to represent
space and time separately. For this purpose we introduce a separate
representation of the column vector, $\bold{U}$, as a column vector
with $(N+1)$ components, each of size $(N+1)^3\times1$: $\bold{U}^T =
[\bold{\hat{U}}_0^T, \ldots,\bold{\hat{U}}_{{N}}^T]$.
Each component $\bold{\hat{U}}_k$, represents the spatial part of the
solution at time point $t_k$.
Let $\bold{\hat{F}}_k$ represent the load vector at time point $t_k$.
In the temporal finite difference approach, at each time step we have
to solve the following linear system,
\begin{equation}
  \label{BE-Spectral}
  \Big(\bold{I_{\text{space}}}+\Delta t \bold{S}_{\text{space}}
  \Big) \bold{\hat{U}}_{k+1} = \bold{\hat{U}}_{k}+ \Delta t
  \bold{\hat{F}}_{k+1}, \quad k=0, \ldots, N.
\end{equation}
  
Here, $\bold{I_\text{space}}$ is the space-identity matrix of size
$(N+1)^3 \times (N+1)^3$, and, $\bold{S}_{\text{space}}$ =
$\bold{A}_{D}$+$\bold{A}_{C}$+$\bold{A}_{R}$, is the space spectral
matrix which is positive definite. Hence, we have,
$\text{det}(\bold{I_\text{space}}+\Delta t \bold{S}_{\text{space}})
\neq 0$, and the linear system Eq. (\ref{BE-Spectral}) has a unique
solution.
We can unite the space and time operators and solve
Eq. \eqref{BE-Spectral} in a single step. First, we rewrite
Eq. \eqref{BE-Spectral} as follows:
\begin{equation}
  \label{TT_BE}
  \frac{1}{\Delta t}\bold{\hat{U}}_{k+1} - \frac{1}{\Delta t}\bold{\hat{U}}_{k} +\bold{S}_{\text{space}} \bold{\hat{U}}_{k+1}= \bold{\hat{F}}_{k+1}.
\end{equation}
Then we define the time derivatives matrix, $\bold{T}$, consistent
with the backward Euler scheme,
\begin{equation}
  \label{Mat_Time}
  \bold{T}=\frac{1}{\Delta t}
  \begin{pmatrix}
    1 & 0 & 0 &\cdots & 0 &0 \\
    - 1 & 1 & 0 & \cdots & 0 & 0\\
    0&  - 1 & 1 & \ldots & 0 & 0\\
    \vdots & \vdots & \vdots & \ddots & \vdots & \vdots\\
    0 & 0 & 0 & \cdots & - 1 &  1
  \end{pmatrix}_{(N+1) \times (N+1)}.
\end{equation}
By employing $\bold{T}$, $\bold{S}_{\text{space}}$, and the
Kronecker product, $\otimes$, we can construct the matricization of
the space-time operator, $\bold{T} \otimes \bold{I}_{\text{space}}
+\bold{I}_t \otimes \bold{S}_{\text{space}}$, where $\bold{I}_{t}$
is the time-identity matrix of size $(N+1) \times (N+1) $.
Hence, the linear system, in the finite differences approach is:
\begin{equation}
  (\bold{T} \otimes \bold{I}_{\text{space}} +\bold{I}_t \otimes \bold{S}_{\text{space}})\bold{U} = \bold{F}. 
\end{equation}
We note that in Ref.~\cite{Dolgov-Khoromskij-Oseledets:2012}, the
authors have employed the same finite difference technique to solve
the heat equation in \TTf{} format.

\PGRAPH{$\bullet$ Time discretization on Chebyshev grid}: When the
order of the Chebyshev polynomials, $N \rightarrow \infty$, the PDE
approximate solution, $u_h(t_n, \bold{x})$, at time $t=t_n$ converges
exponentially in space, as $N^{-|\alpha|}$
\cite[Appendix~A.4]{funaro1997spectral}, and linearly with time step
$\Delta t$,
\begin{equation}
  \label{Order:convergence}
  \|u(t_n,\bold{x}) - u_h(t_n, \bold{x})\|_{L^{2}(\Omega)} \leq (\Delta t + N^{-|\alpha|}) \Big \|\frac{d^{|\alpha|}u(t_n,\bold{x})}{d \bold{x}^{|\alpha|}} \Big\|_{L^2(0,T;L^{2}(\Omega))},
\end{equation}
where $\frac{d^{|\alpha|}u(t,\bold{x})}{d \bold{x}^{|\alpha|}} \in
L^2(\Omega)$, and $u(t, \bold{x}) \in C^0(\Omega)$ for all $t \in
(0,T)$.
Here, $\alpha=(\alpha_1,\alpha_2,\alpha_3)$ is a multi-index, which
characterizes the smoothness of the CDR solution in space by
$|\alpha|:=\alpha_1+\alpha_2+\alpha_3$.
It can be seen that even for higher order Chebyshev polynomials, the
global error is dominated by the error introduced by the temporal
scheme that remains linear, which is the primary drawback of the
temporal finite difference approach.
To recover global exponential convergence, we apply Chebyshev spectral
collocation method for discretization of both space and time
variables.
To accomplish this, we construct a single one-dimensional differential
operator, $\frac{\partial}{\partial t}$ in matrix form, $(\bold{S}_t
)_{ij}$, on temporal Chebyshev grid, with collocation points $t_0,
t_1\ldots, t_{N}$, as follows (see Eq.\eqref{time_derivative}).
\begin{equation}
  (\bold{S}_t)_{ij}=\frac{dl_{j}(t)}{dt} \Big |_{t_i}, \quad 0\leq i, j \leq N.
  \label{time_collocation}
\end{equation}

Hence, the linear system, in the space-time collocation method is,
\begin{equation}
  (\bold{S}_t \otimes \bold{I}_{\text{space}} +\bold{I}_t \otimes \bold{S}_{\text{space}})\bold{U} = \bold{F}, 
  \label{cheb_lin_syst}
\end{equation}
where (see Eq. \ref{mat_eq}) $\mat{A}_t=\bold{S}_t \otimes \bold{I}_{\text{space}}$.
To construct the space-time operator of the CDR equation on the
Chebyshev grid we also need to construct the space part,
$\bold{S}_{\text{space}}$.

\subsubsection{Space discretization on Chebyshev grids}

~\\~\\[-0.5em]
\PGRAPH{$\bullet$ Diffusion Operator on Chebyshev grid}:
\label{Mat:diff}
In this section, we focus on the matricization of the diffusion term,
$ \kappa(t,\bold{x})\Delta $ $\rightarrow$ $\bold{A}_{D}$. The
diffusion operator $\Delta $ is constructed as follows,
\begin{equation}
  \begin{split}
    \label{Spectral_Diff}
    \Delta = \bold{I}_{t}\otimes \bold{S}_{xx} \otimes \bold{I}_{y}\otimes  \bold{I}_{z}+\bold{I}_{t} \otimes \bold{I}_{x} \otimes \bold{S}_{yy}  \otimes\bold{I}_{z}  + \bold{I}_{t} \otimes \bold{I}_{x} \otimes  \bold{I}_{y} \otimes\bold{S}_{zz}.
  \end{split}
\end{equation}
Then, the function $\kappa(t,\bold{x})$ is incorporated to form the
diffusion term $\bold{A}_{D}$
\begin{equation}
  \label{eqn:AD operator}
  \bold{A}_{D} = diag(\boldsymbol{K})\Delta,
\end{equation}
where $diag({\ldots})$ denotes a diagonal matrix, and
${\boldsymbol{K}}$ is a vector of size $(N+1)^4$ containing the
evaluation of $\kappa(t,\bold{x})$ on the Chebyshev space-time grid.

\PGRAPH{$\bullet$ Discretization of the convection term on Chebyshev
  grid}:
Here, we focus on the matricization of the convection term,
$\vec{b}(t,\bold{x})\cdot\nabla\rightarrow \bold{A}_C$, with the
convective function, $\vec{b}(t,x,y,z)$, which we assume in the form
\begin{equation}
  \vec{b}(t,x,y,z)=[b^x(t,x,y,z)~b^y(t,x,y,z)~b^z(t,x,y,z)].
\end{equation}
Then, we construct the convection term $\bold{A}_C$
\begin{equation}
  \label{convec:mat}
  \begin{split}
    \bold{A}_{C}= \ & diag(\bold{B}^x) \left(\bold{I}_{t} \otimes \bold{S}_{x} \otimes \bold{I}_{y}\otimes\bold{I}_{z}\right) + \ldots\\
    & diag(\bold{B}^y)\left(\bold{I}_{t} \otimes \bold{I}_{x}\otimes \bold{S}_{y}  \otimes\bold{I}_{z}\right) + \ldots \\
    & diag(\bold{B}^z)\left(\bold{I}_{t} \otimes \bold{I}_{x} \otimes  \bold{I}_{y}\otimes \bold{S}_{z}\right),
  \end{split}
\end{equation}
where $\bold{B}^x,\bold{B}^y$, and $\bold{B}^z$ are vectors of size
$(N+1)^4$ containing the evaluation of the functions $b^x,\ b^y$ and
$b^z$ on the Chebyshev space-time grid.

\PGRAPH{$\bullet$ Discretization of the reaction term on Chebyshev
  grid}: Here, we focus on the matricization of the reaction term,
$c(t,\bold{x})$ $\rightarrow$ $\bold{A}_{R}$, which is given by
$\bold{A}_{R} = diag(\bold{C})\bold{I}_{t} \otimes \bold{I}_{x}
\otimes \bold{I}_{y}\otimes\bold{I}_{z}$, where $\bold{C}$ is the
evaluation of the function $c(t,\bold{x})$ on the Chebyshev space-time
grid.

\subsubsection{Initial and boundary conditions on space-time Chebyshev grids} 
\label{Bd_cond}

\begin{figure}[!t]
  \centering
  \includegraphics[width=0.5\textwidth]{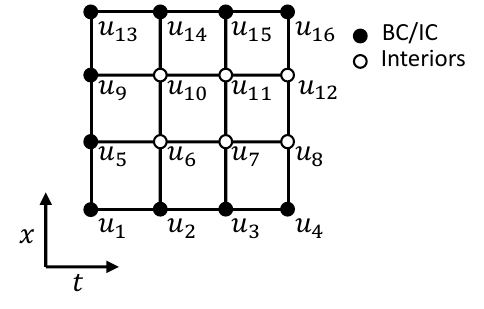}
  \caption{1D Space time grid with $N=4$ collocation nodes.}
  \label{fig:ST-2d-grid}
\end{figure}

So far, we have not incorporated the boundary conditions (BC) and the
initial condition (IC), given in \eqref{model:prob:bdy}, into the
linear system $\bold{A}\bold{U} = \bold{F}$ \eqref{cheb_lin_syst}.
In the space-time method, we consider the IC equivalent to the BC.
The nodes of the Chebyshev grid are split into two parts: $(i)$ BC/IC
nodes, and $(ii)$ interior nodes, where the solution is unknown.
Let $i^{\text{Bd}}$, and $i^{\text{Int}}$ be the set of multi-indices,
introduced in \cite{dolgov2021guaranteed},
the BC/IC and interior nodes, respectively.

Boundary and initial conditions are imposed by explicitly enforcing
the BC/IC nodes to be equal to the BC, $g(t,x)$, or to the IC, $h(x)$,
and then reducing the linear system for all nodes into a smaller
system only for the interior nodes,
\begin{equation}
  \label{eqn:reduced_system}
  \bold{A}^{\text{Int}} \bold{U}^{\text{Int}}=\bold{F}^{\text{Int}} -\bold{F}^{\text{Bd}}.
\end{equation}

To make it clear we consider below a simple example with $N = 4$
collocation points in two dimensions $(t,x)$.
The set of Chebyshev nodes can be denoted by multi-indices as,
$\bold{U}:=\{u_1,u_2,\ldots,u_{16} \} $ (Figure
\ref{fig:ST-2d-grid}). After imposing the BC/IC conditions, the linear
system with BC/IC nodes, $\bold{A}\bold{U} = \bold{F}$, is

\begin{equation}
  \label{Full_system}
  \begin{split}
    &u_1 = g (t_{0},{x}_{0})\\
    &\quad \vdots \\
    &u_4 = g (t_{3},{x}_{0})\\
    &u_5 = h(x_1), \\
    &A_{6,1} u_1 + A_{6,2} u_2+ \ldots +A_{6,15}u_{15} + A_{6,16}u_{16} =F_6\\
    &\quad \vdots \\
    &A_{8,1} u_1 + A_{8,2} u_2+ \ldots +A_{8,15}u_{15} + A_{8,16}u_{16} =F_8\\
    &u_9 = h(x_{2})\\
    &A_{10,1} u_1 + A_{10,2} u_2+ \ldots +A_{10,15}u_{15} + A_{10,16}u_{16} =F_{10}\\
    &\quad \vdots \\
    &A_{12,1} u_1 + A_{12,2} u_2+ \ldots +A_{12,15}u_{15} + A_{12,16}u_{16} =F_{12}\\
    & \quad \vdots \\
    &u_{13} = g (t_{0},{x}_{3})\\
    &\quad \vdots \\
    &u_{16} = g (t_{3},{x}_{3}).\\
  \end{split}
\end{equation}

The unknown values $(u_6,u_7,u_8,u_{10},u_{11},u_{12})$ associated to
six interior nodes, satisfy the reduced linear system whose equations
read as
\begin{multline}
  \qquad
  A_{l,6} u_6+A_{l,7}u_7+A_{l,8}u_{8} + A_{l,10}u_{10} +A_{l,11}u_{11} +A_{l,12}u_{12}\\
  =F_l-A_{l,1}g (t_0,x_0)-A_{l,2}g (t_1,x_0)-\ldots,
  \qquad
  \label{eqn:sub_system_BC}
\end{multline}
with $l \in \{6,7,8,10,11,12\}$, where we transfer the values of BC/IC
nodes to the right-hand side.

\section{Tensor networks}
\label{sec:Tensor-Networks}

In this section, we introduce the \TTf{} format, the specific tensor
network we use in this work, the representation of linear operators in
TT-matrix format, and the cross-interpolation method.
All these methods are fundamental in the tensorization of our spectral
collocation discretization of the CDR equation.
For a more comprehensive understanding of notation and concepts, we
refer the reader to the following references:
Refs.~\cite{kolda2009tensor,oseledets2011tensor,dolgov2021guaranteed}, which
provide detailed explanations.

\subsection{Tensor train}
\label{subsec:TNs}
\begin{figure}[h] 
  \centering
  \includegraphics[width=0.75\textwidth]{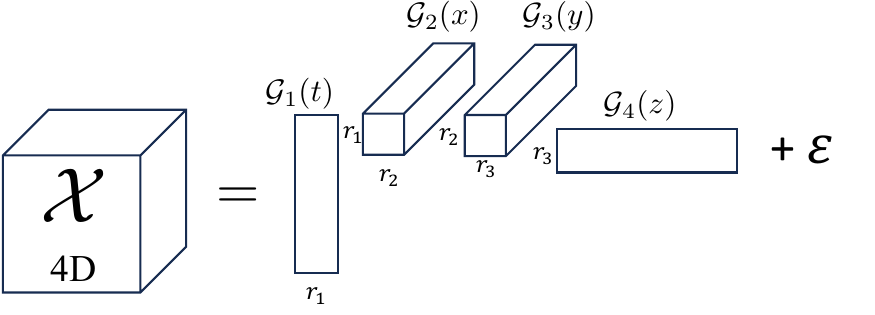}
  \caption{ \TTf{} decomposition of a 4D tensor $\ten{X}$, with \TTf{} rank,
    $\mathbf{r} = \big[r_1,r_2,r_3\big]$, and approximation error
    $\varepsilon$, in accordance with Eq.~\eqref{eqn:TT_def_element}.}
   \label{fig:TT_4D}
\end{figure}

The TT format, introduced by Oseledets in
2011~\cite{oseledets2011tensor}, represents a sequential chain of
matrix products involving both two-dimensional matrices and
three-dimensional tensors, referred to as TT-cores.
We can visualize this chain as in Figure \ref{fig:TT_4D}.
Given that tensors in our formulation are at most four dimensional
(one temporal and three spatial dimensions), we consider the tensor
train format in the context of 4D tensors.
Specifically, the \TTf{} approximation $\ten{X}^{TT}$ of a
four-dimensional tensor $\ten{X}$ is a tensor with elements
\begin{align}
  \ten{X}^{TT}(i_1,i_2,i_3,i_4)
  =  \sum^{r_1}_{\alpha_{1}=1}\sum^{r_2}_{\alpha_{2}=1}\sum^{r_3}_{\alpha_{3}=1}
  \ten{G}_1(1,i_1,\alpha_1)\ten{G}_2(\alpha_1,i_2,\alpha_2)\ten{G}_3(\alpha_2,i_3,\alpha_3)\ten{G}_4(\alpha_{3},i_4,1)
  + \varepsilon,
  \label{eqn:TT_def_element}
\end{align}
where the error, $\varepsilon$, is a tensor with the same dimensions
as $\ten{X}$, and the elements of the array $\mat{r} = [r_1,r_2,r_3]$
are the TT-ranks, and quantify the compression effectiveness of the TT
approximation.
Since, each \TTf{} core, $\ten{G}_p(i_k)$, only depends on a single index
of the full tensor $\ten{X}$, e.g., $i_k$, the \TTf{} format effectively
embodies a discrete separation of variables \cite{bachmayr2016tensor}.
In Figure~\ref{fig:TT_4D} we show a four-dimensional array
$\ten{X}(t,x,y,z)$, decomposed in TT-format.

\subsection{Linear Operators in TT-matrix format}
\label{SUB:Lin_operinTT}
Suppose that the approximate solution of the CDR equation is a 4D
tensor $\ten{U}$, then the linear operator $\ten{A}$ acting on that
solution is represented as an 8D tensor.
The transformation $\ten{A}\ten{U}$ is defined as:
\begin{align*}
  \big(\ten{A}\ten{U}\big)(i_1,i_2,i_3,i_4)
  = \sum_{j_1,j_2,j_3,j_4} \ten{A}(i_1,j_1,\dots,i_4,j_4)\ten{U}(j_1,\dots,j_4).
\end{align*}
The tensor $\ten{A}$ and the matrix operator $\bold{A}$, defined in
Eq.~\eqref{eqn:reduced_system}, are related as:
\begin{equation}
  \ten{A}(i_1,j_1,\dots,i_4,j_4) = \bold{A}(i_1i_2i_3i_4,j_1j_2j_3j_4).
  \label{eqn:ten_mat_op_relation}
\end{equation}
Thus, we can construct the tensor $\ten{A}$ by suitably reshaping and
permuting the dimensions of the matrix $\bold{A}$.
The linear operator, $\ten{A}$, can be further represented in a
variant of \TTf{} format, called \emph{TT-matrix},
cf.~\cite{truong2023tensor}.
The component-wise \TTf{}-matrix $\ten{A}^{TT}$ is defined as:
\begin{equation}
  \ten{A}^{TT}(i_1,j_1,\ldots,i_4,j_4)
  =  \sum_{\alpha_{1},\alpha_2,\alpha_{3}}
  \ten{G}_1\big(1,(i_1,j_1),\alpha_1\big)
   \ldots\ten{G}_4\big(\alpha_{3},(i_4,j_4),1\big),
  \label{eqn:TT-matrix-componentwise}
\end{equation}
where $\ten{G}_{k}$ are 4D TT-cores .
\begin{figure}[h]
  \centering
  \includegraphics[width=0.8\textwidth]{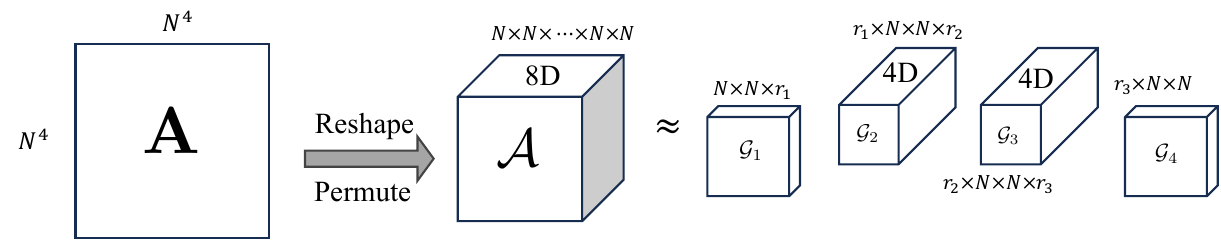}
  \caption{Representation of a linear matrix $\mat{A}$ in the
    TT-matrix format. First, we reshape the operation matrix
    $\bold{A}$ and permute its indices to create the tensor
    $\ten{A}$. Then, we factorize the tensor in the tensor-train
    matrix format according to Eq.~\eqref{eqn:TT-matrix-componentwise}
    to obtain $\ten{A}^{\TTf{}}$.}
  \label{fig:TT-matrix-format}
\end{figure}
Figure~\ref{fig:TT-matrix-format} shows the process of transforming a
matrix operator $\bold{A}$ to its tensor format, $\ten{A}$, and,
finally, to its TT-matrix format, $\ten{A}^{\TTf{}}$.

We can further simplify the TT-matrix representations of the matrix
$\bold{A}$ if it is a Kronecker product of matrices, i.e.
$\bold{A}=\bold{A}_1\otimes\bold{A}_2\otimes\bold{A}_3\otimes\bold{A}_4$.
Based on the relationship defined in
Eq.~\eqref{eqn:ten_mat_op_relation}, the tensor $\ten{A}$ can be
constructed using tensor product as $\ten{A} = \bold{A}_1 \tenprod
\bold{A}_2 \tenprod \bold{A}_3 \tenprod \bold{A}_4$.
This implies the internal ranks of the TT-format of $\ten{A}$
in \eqref{eqn:TT-matrix-componentwise} are all equal to $1$.
In such a case, all summations in
Eq.~\eqref{eqn:TT-matrix-componentwise} reduce to a sequence of single
matrix-matrix multiplications, and the \TTf{} format of $\ten{A}$ becomes
the tensor product of $d$ matrices:
\begin{align}
  \ten{A}^{TT} &= \mat{A}_1 \tenprod \mat{A}_2 \tenprod \dots \tenprod \mat{A}_d.
  \label{eqn:TT_matrices}
\end{align}
This specific structure appears quite often in the matrix
discretization, and will be exploited in the tensorization to
construct efficient \TTf{} format.

\subsection{TT Cross Interpolation}
\label{sec:Cross Interpolation}
\begin{figure}[h]
  \centering
  \includegraphics[width=0.5\textwidth]{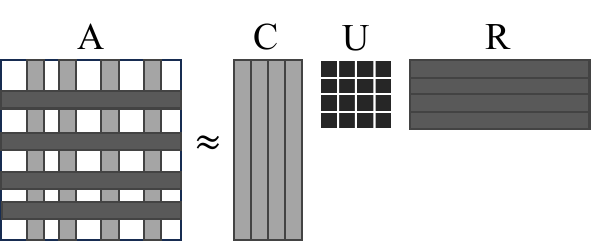}
  \caption{CUR decomposition}
  \label{fig:CUR}
\end{figure}

The original \TTf{} algorithm is based on consecutive applications of
singular value decompositions (SVDs) on unfoldings of a tensor
\cite{oseledets2011tensor}.
Although known for its efficiency, \TTf{} algorithm requires access to the
full tensor, which is impractical and even impossible for extra-large
tensors.
To address this challenge, the cross interpolation algorithm,
TT-cross, has been developed \cite{oseledets2010tt}.
The idea behind TT-cross is essentially to replace the SVD in the TT
algorithm with an \emph{approximate} version of the skeleton/CUR
decomposition \cite{goreinov1997theory, mahoney2009cur}.
CUR decomposition approximates a matrix by selecting a few of its
columns $\mat{C}$, a few of its rows $\mat{R}$, and a matrix $\mat{U}$
that connects them, as shown in Fig.~\ref{fig:CUR}.
Mathematically, CUR decomposition finds an approximation for a matrix
$\mat{A}$, as $\mat{A}\approx \mat{C}\mat{U}\mat{R}$.
The TT-cross algorithm, utilizes the Maximum Volume Principle
(\emph{maxvol algorithm})
\cite{goreinov2010find,mikhalev2018rectangular} to determine
$\mat{U}$.
The maxvol algorithm chooses a few columns, $\mat{C}$ and rows,
$\mat{R}$, of $\mat{A}$, such that, the intersection matrix
$\mat{U}^{-1}$ has maximum volume \cite{savostyanov2011fast}.

TT-cross interpolation and its versions can be seen as an heuristic
generalization of CUR to tensors
\cite{oseledets2008tucker,sozykin2022ttopt}.
TT-cross utilizes the maximum volume algorithm iteratively, often
beginning with few randomly chosen fibers, to select an optimal number
of specific tensor fibers that capture essential information of the
tensor \cite{savostyanov2014quasioptimality}.
These fibers are used to construct a lower-rank \TTf{}
representation. The naive generalization of CUR is proven to be
expensive, which led to the development of various heuristic
optimization techniques, such as, TT-ALS \cite{holtz2012alternating},
DMRG \cite{oseledets2012solution,savostyanov2011fast}, and AMEN
\cite{dolgov2014alternating}.

To solve the CDR equation, we use TT-cross to build the \TTf{} format,
directly from the coefficient functions $\kappa(t,\bold{x})$,
$\bold{b}(t,\bold{x})$, $c(t,\bold{x})$, boundary conditions, initial
conditions, and loading functions.

\subsection{Tensorization}
\label{tensorization}
The space-time discretization produces a linear system for all
interior nodes as specified in Eq.~\eqref{eqn:reduced_system}.
Here, we simplify the notation, and refer to this equation as
$\bold{A}\bold{U}=\bold{F} -\bold{F}^{\text{Bd}}$.

Tensorization is a process of building the \TTf{} format of all components
of this linear system
\begin{equation}
  \bold{A}\bold{U}=\bold{F} -\bold{F}^{\text{Bd}} \rightarrow \ten{A}^{TT}\ten{U}^{TT} = \ten{F}^{TT} - \ten{F}^{\text{Bd},TT},
\end{equation}
where $\ten{A}^{TT} = \ten{S}_t^{TT} + \ten{A}_D^{TT} + \ten{A}_C^{TT}+ \ten{A}_R^{TT}$.

In the matrix form, the operator $\bold{A} = \bold{A}_t + \bold{A}_D +
\bold{A}_C + \bold{A}_R$ is designed to act on the vectorized solution
$\bold{U}$.
In the full tensor format, the solution is kept in its original tensor
form $\ten{U}$, which is a 4D tensor.
Consequently, all the operators $\ten{A}_t,\ \ten{A}_D, \ \ten{A}_C,
\text{ and } \ten{A}_R$ are 8D tensors.
Lastly, $\ten{F}$ and $\ten{F}^{\text{Bd}}$ are both 4D tensors.
Here, given that these operators in their matrix form have Kronecker
product structure, their \TTf{} format can be constructed by using
component matrices as \TTf{} cores.
To construct the tensor format of the operators acting on the interior
nodes, we need to define some sets of indices:
\begin{equation*}
  \begin{split}
    & \ten{I}_t = 1:N \ \text{index set for time variable},\\
    & \ten{I}_s = 1:(N-1) \ \text{index set for space variable}.
  \end{split}
\end{equation*}

\PGRAPH{$\bullet$ TT-matrix Time Operator}, $\ten{A}^{TT}_t$: The temporal operator
in TT-matrix format acting only on the interior nodes is constructed
as:
\begin{equation}
  \ten{A}^{TT}_t = \bold{S}_t(\ten{I}_t,\ten{I}_t)\tenprod \bold{I}_{N-1} \tenprod \bold{I}_{N-1} \tenprod \bold{I}_{N-1},
  \label{eqn:At_TT_op}
\end{equation}
where $\bold{I}_{N-1}$ is the identity matrix of size
$(N-1)\times(N-1)$.
  
\PGRAPH{$\bullet$ TT-matrix Diffusion Operator}, $\ten{A}^{TT}_D$: The Laplace
operator in TT-matrix format is constructed as:
\begin{equation}
  \begin{split}
      \Delta^{TT}=\bold {I}_{N}\tenprod \bold{S}_{xx}(\ten{I}_s,\ten{I}_s) \tenprod \bold{I}_{N-1}\tenprod \bold{I}_{N-1}
      + &\bold{I}_{N} \tenprod \bold{I}_{N-1} \tenprod \bold{S}_{yy}(\ten{I}_s,\ten{I}_s)  \tenprod\bold{I}_{N-1}\\
      + &\bold{I}_N \tenprod \bold{I}_{N-1} \tenprod  \bold{I}_{N-1} \tenprod\bold{S}_{zz}(\ten{I}_s,\ten{I}_s).
  \end{split}
\end{equation}

Hence, the Laplace operator is a sum of three linear operators in
TT-matrix formats, which again is a TT-matrix.
Further, we need to format the diffusion coefficient
$\kappa(t,\bold{x})$ on the Chebyshev grid of interior nodes in TT
format.
To transform $\kappa(t, \textbf{x})\rightarrow \ten{K}^{TT}$, we
applied the cross-interpolation technique described in
Section~\ref{sec:Cross Interpolation}.
Finally, we need to transform $\ten{K}^{TT}$ in TT-matrix format to
be able to multiply it with $\Delta^{TT}$.
The transformation to \TTf{} $\rightarrow$ TT-matrix is given in
Algorithm \ref{alg:Ktt_to_Koptt}.

\begin{algorithm}
  \caption{Reformat TT, $\tentt{K}$, in TT-matrix format, $\ten{K}^{op,TT}$ }
  \label{alg:Ktt_to_Koptt}
  \KwData{$\tentt{K}$}
  \KwResult{$\ten{K}^{op,TT}$}
  Get the cores $\left\{\ten{G}_k\right\}_{k=1,2,3,4}$ from $\tentt{K}$.\\
  $\{r_0=1,r_1,r_2,r_3,r_4=1\}$ are \TTf{} ranks.\\
  $\{n_1,n_2,n_3,n_4\}$ are tensor sizes.\\
  \For{k = 1:4}{
    $\ten{G}new_k = zeros(r_{k-1},n_k,n_k,r_k)$ \tcp{create a zero tensor core}
    \For{i=1:$r_{k-1}$}{
      \For{j=1:$r_k$}{
        $\ten{G}new_k(i,:,:,j) = diag(\ten{G}_k(i,:,j))$
      }
    }
  }
  Build $\ten{K}^{op,TT}$ using $\{\ten{G}new_k\}$ as its cores.
\end{algorithm}
  
\PGRAPH{$\bullet$ TT-matrix Convection Operator}, $\ten{A}^{TT}_C$:
The convection operator in TT-matrix format is constructed as:
\begin{equation}
  \begin{split}
    \nabla^{TT}_x = \bold {I}^{N}\tenprod \bold{S}_{x}(\ten{I}_s,\ten{I}_s) \tenprod \bold{I}^{N-1}\tenprod \bold{I}^{N-1},\\[0.5em]
    \nabla^{TT}_y = \bold {I}^{N}\tenprod \bold{I}^{N-1}\tenprod\bold{S}_{y}(\ten{I}_s,\ten{I}_s) \tenprod \bold{I}^{N-1}, \\[0.5em]
    \nabla^{TT}_z = \bold {I}^{N} \tenprod \bold{I}^{N-1}\tenprod \bold{I}^{N-1} \tenprod\bold{S}_{z}(\ten{I}_s,\ten{I}_s).\\
  \end{split}
\end{equation}

Next, the tensor operators $\ten{B}^{op,TT}_x$, $\ten{B}^{op,TT}_y$
and $\ten{B}^{op,TT}_z$ are computed from the functions $b^x,\ b^y$
and $b^z$ in the same way as being computed with $\kappa(t,\bold{x})$.
Then the TT-matrix format of convection operator is constructed as:
\begin{equation}
  \ten{A}_C^{TT} =
  \ten{B}_x^{op,TT}\nabla_x^{TT} +
  \ten{B}_y^{op,TT}\nabla_y^{TT} +
  \ten{B}_z^{op,TT}\nabla_z^{TT}.
\end{equation}

\PGRAPH{$\bullet$ TT-matrix Reaction Operator}, $\ten{A}^{TT}_R$: The
reaction operator $\ten{A}^{TT}_R$ basically is the function
$c(t,\bold{x})$ being converted to an operator in the same way
as with other coefficient functions.

Then the TT-matrix format of the operator $\ten{A}$ is constructed as
\[\ten{A}^{TT} = \ten{A}^{TT}_t + \ten{A}^{TT}_D + \ten{A}^{TT}_C + \ten{A}^{TT}_R.\]

\PGRAPH{$\bullet$ TT-matrix Loading Tensor}, $\ten{F}^{TT}$: The TT-matrix
loading tensor $\ten{F}^{TT}$ is constructed by applying the cross
interpolation on the function $f(t,\bold{x})$ on the grid of interior
nodes.
  
\PGRAPH{$\bullet$ TT Boundary Tensor}, $\ten{F}^{Bd,TT}$: The boundary
tensor $\ten{F}^{Bd}$ is used to incorporate the information
about boundary and initial conditions into the linear system of
the interior nodes.
Only in the case of homogeneous BC and zero IC, this tensor is a
zero tensor.
Basically, the idea is to construct an operator $\ten{A}^{map}$,
similar to $\ten{A}$.
The difference is that $\ten{A}^{map}$ will be a mapping from all
nodes to only unknown nodes.
The size of the operator $\ten{A}^{map}$ is $N\times (N+1) \times
(N-1)\times (N+1) \times (N-1)\times (N+1) \times(N-1)\times
(N+1)$.
Then, the boundary tensor $\ten{F}^{Bd}$ is computed by apply
$\ten{A}^{map}$ to a tensor $\ten{G}^{Bd}$ containing only the
information from BC/IC.
The details about constructing $\ten{A}^{map,TT}$ and
$\ten{G}^{Bd,TT}$ are included in the Appendix \ref{APP:A_map_tt}.

At this point, we have completed constructing the TT-format of the
linear system $\ten{A}^{TT}\ten{U}^{TT} = \ten{F}^{TT} -
\ten{F}^{Bd,TT}$.
To solve the \TTf{} linear system by optimization techniques, we used
the routines \texttt{amen\_cross, amen\_solve} and \texttt{amen\_mm}
form the MATLAB TT-Toolbox \cite{Oseledets:2023}.


\section{Numerical Results}
\label{results}

We investigate the computational efficiency and potential low-rank
structures of the space-time solver by implementing the space-time
operator in the tensor train format.
To validate the algorithm presented earlier and assess its
performance, we conducted several numerical experiments using MATLAB
R2022b on a 2019 Macbook Pro 16 equipped with a 2.4 GHz 8-core i9 CPU.
All our implementations rely on the MATLAB TT-Toolbox
\cite{Oseledets:2023}.
In the first example, we consider a manufactured solution of the
time-dependent CDR equation with constant coefficients.
The second example involves the space-time solution operator with
non-constant coefficients. Finally, the third example showcases the
behavior of the methods applied to the CDR with a non-smooth
manufactured solution.

For all the examples, we compare the \TTf{} format of the space-time
operator with the Finite Difference-Finite Difference (\FDFD) method
and the Spectral-Spectral (\SPSP) collocation method.
Additionally, we compared the \TTf{} format of \SPSP{} with the full grid
formulation of the \SPSP{} method. To quantify the compression achieved,
we introduced the compression ratio \CmpRate, defined as:
\begin{equation}
  \label{eq:CR}
  \CmpRate:=\frac{\Mem(\ten{X}^{\TTf})}{\Mem(\ten{X})},
\end{equation}
where $\Mem(\ten{X}^{TT})$ and $\Mem(\ten{X})$ are the total memory
required for storing the scheme's unknows in the \TTf{} format, e.g.,
$\ten{X}^{\TTf}$ and the full tensor format, e.g., $\ten{X}$.
We will show that the \TTf{} format may achieve significant
compression of spectral-collocation differential operators from
terabyte sizes to megabytes and even kilobytes when solving the CDR
equation on very fine meshes.
Such a compression especially reduces the memory requirements for
solving the linear system that results from the CDR equation
discretization.

\subsection{Test~1: Manufactured Solution with Constant Coefficients}

In this first test, we study the convergence behavior of the benchmark
problem described by Eq. \eqref{eq:model:problem} on the computational
domain $[0,1]\times[-50,50]^3$.
We consider constant coefficients, specifically
$\kappa(t,\bold{x})=1$, $\bold{b}(t,\bold{x})=(1,1,1)$, and
$c(t,\bold{x})=1$.
The right-hand side force function $f(t,\bold{x})$ and the boundary
conditions are determined in accordance with the exact solution
$u(t,x,y,z):=\sin(2 \pi (t+x+y+z))$.

\begin{figure}[ht]
  \centering
  \includegraphics[width = 0.85\textwidth]{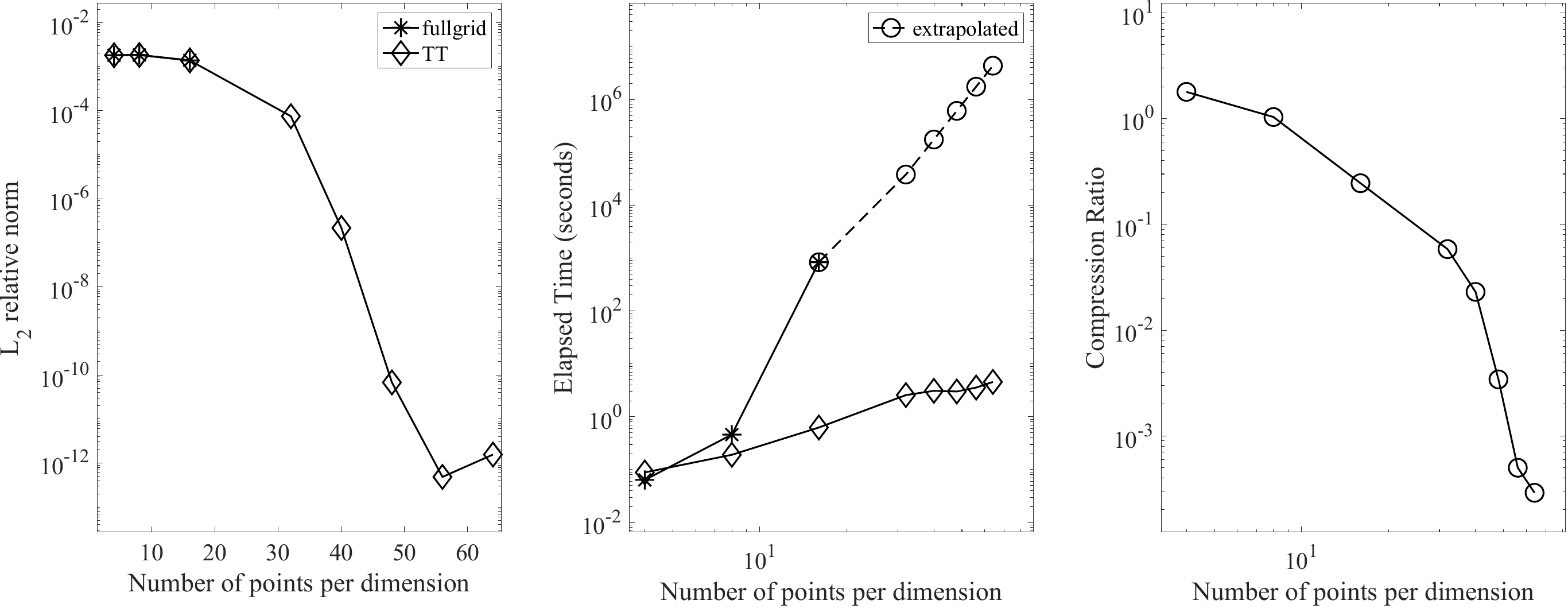}
  \caption{Test 1: Left Panel: Relative error curve in $L^2$ norm
    showing the exponential convergence of \SPSP{} schemes. Middle
    Panel: Elapsed time in seconds.  Right Panel: Compression ratio of
    the solution. All plots are versus the number of points per
    dimension.}
  \label{fig:ConstantCDR_SPTT_SPFG}
\end{figure}

Figure~\ref{fig:ConstantCDR_SPTT_SPFG} presents the results of this
test.
The left panel displays the relative error curve in the $L^2$ norm,
demonstrating the exponential convergence of the \SPSP{} schemes.
The middle panel depicts the elapsed time in seconds, indicating that
the \TTf{} format enables solutions with much higher resolution
compared to the full grid format. The right panel illustrates the
compression ratio of the solution. We plot all curves against the
number of points per dimension.

In Figure~\ref{fig:ConstantCDR_SPTT_SPFG}, we compare the \SPSP{} method
in the full grid format and the \SPSP{} method in \TTf{} format.
In the left panel, we plot the $L^2$ relative norms computed using the
full grid and \TTf{} format against the number of nodes per dimension.
The \TTf{} format provides the same level of accuracy as the full grid
scheme.
However, the \TTf{} format allows us to handle large data sets, unlike
the full grid approach, which is limited to small grid data.
This characteristic ensures that the \TTf{} format utilizes
significantly less memory than the full-grid approach, making handling
large data sets practical.
Furthermore, the numerical results demonstrate the well-known
exponential convergence of the spectral method.
In the full grid computation, we achieve an accuracy of approximately
$10^{-3}$, while in the \TTf{} format, we attain an accuracy of
approximately $10^{-12}$ with only a minimal increase in grid data.
At 64 nodes, the accuracy does not improve further as it reaches the
truncation tolerance of the \TTf{} format.

In the middle panel, we compare the elapsed time (in seconds) required
for both the \TTf{} format and full grid approaches.
In the \TTf{} format, the elapsed time to compute the space-time
operator increases linearly while preserving the invariance of the
variables.
On the other hand, in the full grid approach, the elapsed time
increases exponentially.
This behavior aligns with the theoretical prediction as the complexity
grows linearly.
For instance, with 16 nodes per dimension, the time taken in the
\TTf{} format (represented by the third cross point) is 1400 times
faster than in the full grid format (represented by the circled point).

Finally, we emphasize that this speed-up will further improve for
larger grids.
For example, based on the extrapolated elapsed time for the full grid
scheme, we expect the \TTf{} method to be approximately $9\times10^5$
times faster with 64 nodes per dimension.

Another interesting aspect is the compression ratio, which we defined
in Eq.~\eqref{eq:CR}.
In the right panel of Figure~\ref{fig:ConstantCDR_SPTT_SPFG}, we
present the plot of the compression ratio versus the number of nodes
per dimension, indicated by the curve.
The compression ratio exponentially decreases as the number of nodes
per dimension increases.
At $64$ nodes per dimension, we achieve a compression ratio that saves
approximately $10^{4}$ orders of magnitude in storage.
This remarkable reduction in storage requirements is further
illustrated in Table \ref{tab:test_1_compress}, where we showcase the
compression of the aggregated operator $\ten{A}$, which is the memory
bottleneck of the full grid scheme.
The \TTf{} format allows full utilization of terabyte-sized operators
while utilizing only kilobytes of storage.
These results demonstrate the computational advantages of the \SPSP{}
scheme in \TTf{} format compared to existing techniques.

\begin{table}[h]
  \centering
  \begin{tabular}{ |c|c|c|c|} 
    \hline Number of points per dimension & size of $\ten{A}$ in TB &
    size of $\ten{A}^{TT}$ in KB & Compress ratio \\ \hline 8&1.66e-5
    &3.58e0 & 2.00e-4\\ 16&1.23e-2 &1.88e1 & 1.42e-6\\ 32&5.10e0
    &8.53e1 & 1.56e-8\\ 64&1.64e3 &3.63e2 & 2.06e-10\\ \hline
  \end{tabular}
  \caption{Storage cost comparison of $\ten{A}$ between full-format
    and \TTf-format formulations. At $64$ points, the \TTf-format
    allows a compression from $1640$ Tb to $363$ Kb, which is about
    ten orders of magnitude.}
  \label{tab:test_1_compress}
\end{table}

\begin{figure}[ht]
  \centering
  \includegraphics[width=0.9\textwidth]{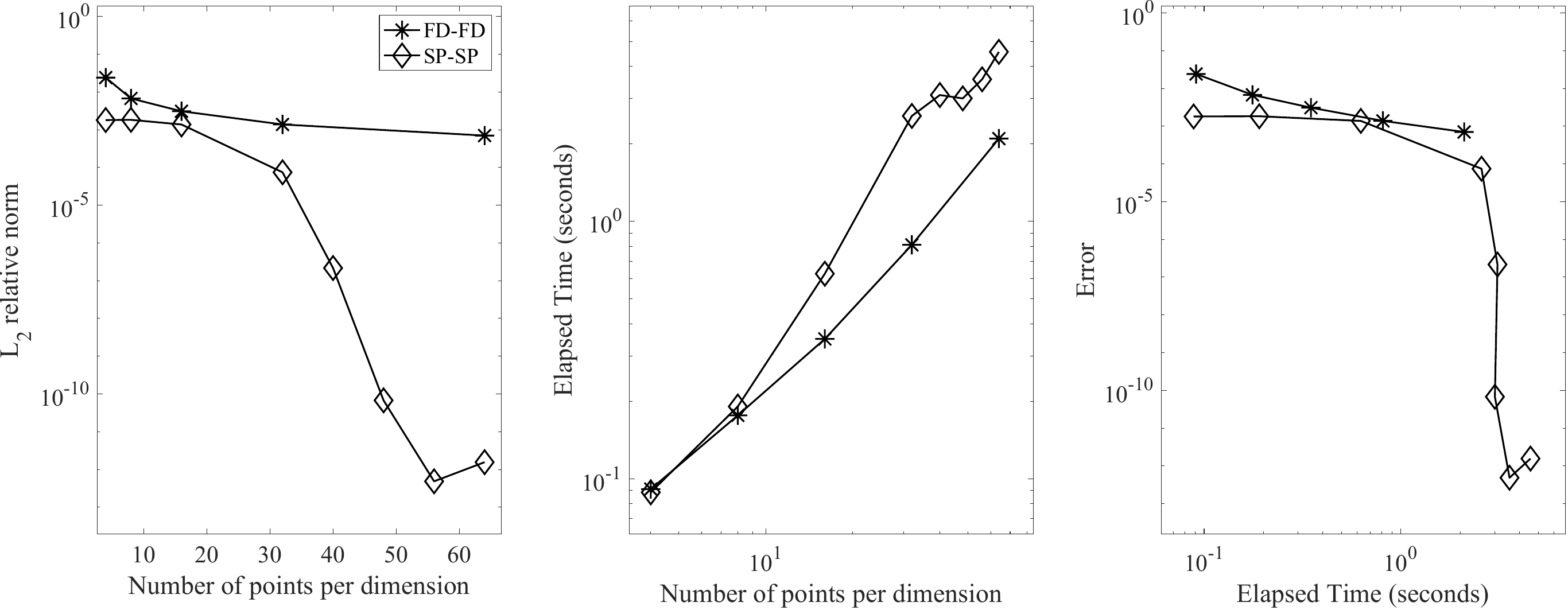}
  \caption{Test 1: Comparison between \TTf{} \SPSP{} and \TTf{} \FDFD.
    Left Panel: Plot of relative error versus number of points per
    dimension. Middle Panel: Plot of elapsed time versus number of
    points per dimension. Right Panel: Plot of relative error versus
    elapsed times, which shows that the \TTf{} \SPSP{} is more
    efficient compared to \TTf{} \FDFD.}
  \label{ConstantCDR_SPTT_FDTT}
\end{figure}

Next, we aim to compare the performance of the \SPSP{} scheme with the
\FDFD{} scheme, both in \TTf{} format.
The space-time \FDFD{} scheme was initially proposed by Dolgov et
al. \cite{Dolgov-Khoromskij-Oseledets:2012}.
In Fig.~\ref{ConstantCDR_SPTT_FDTT}, we compare the \SPSP{} scheme
in \TTf{} format and the \FDFD{} scheme in \TTf{} format.
The left panel shows the $L^2$ relative error norms computed by both schemes in TT format.
As the \SPSP{} scheme converges exponentially, it provides
significantly more accurate results than the \FDFD{} scheme, which
converges linearly.
For instance, at $64$ nodes per dimension, the approximation error of
the \SPSP{} scheme is approximately $10^{-11}$, while the
approximation error of the \FDFD{} scheme is only of the order of
$10^{-3}$.
This fact shows that our \TTf{}-based spectral-collocation method is
approximately $10^{-8}$ times more accurate than the \FDFD{}-based for
the same grid data.

In the middle panel, we plot the elapsed time (in seconds) against the
number of nodes per dimension.
At $64$ nodes per dimension, the \SPSP{} scheme takes roughly twice as
long as the \FDFD{} scheme.
This indicates that the \SPSP{} scheme is slower by a factor of two
compared to the \FDFD{} scheme.
However, in the right panel, we observe that the \SPSP{} scheme
computes highly accurate solutions relative to the consumed time, in
contrast to the \FDFD{} scheme.
Specifically, the \SPSP{} scheme requires $10$ seconds to compute a
solution accurate to the order of $10^{-12}$, whereas the \FDFD{}
scheme computes a solution accurate to the order of $10^{-3}$ within
the same time frame.

In conclusion, we find that the \SPSP{} scheme computes highly
accurate approximate solutions compared to existing techniques, such
as the \FDFD{} scheme.

\subsection{Test~2: Manufactured Solution with variable coefficients}

In this test case, we examine the same model problem as above with the
variable coefficients $\kappa(t,\bold{x})=\exp(-t^2)$,
$\bold{b}(t,\bold{x})=(\sin(2 \pi x), \cos(2 \pi y), \sin(2 \pi z))$,
and $c(t,\bold{x})=\cos(2 \pi (t+x+y+z))$.
\begin{figure}[ht]
  \centering
  \includegraphics[width = 0.9 \textwidth]{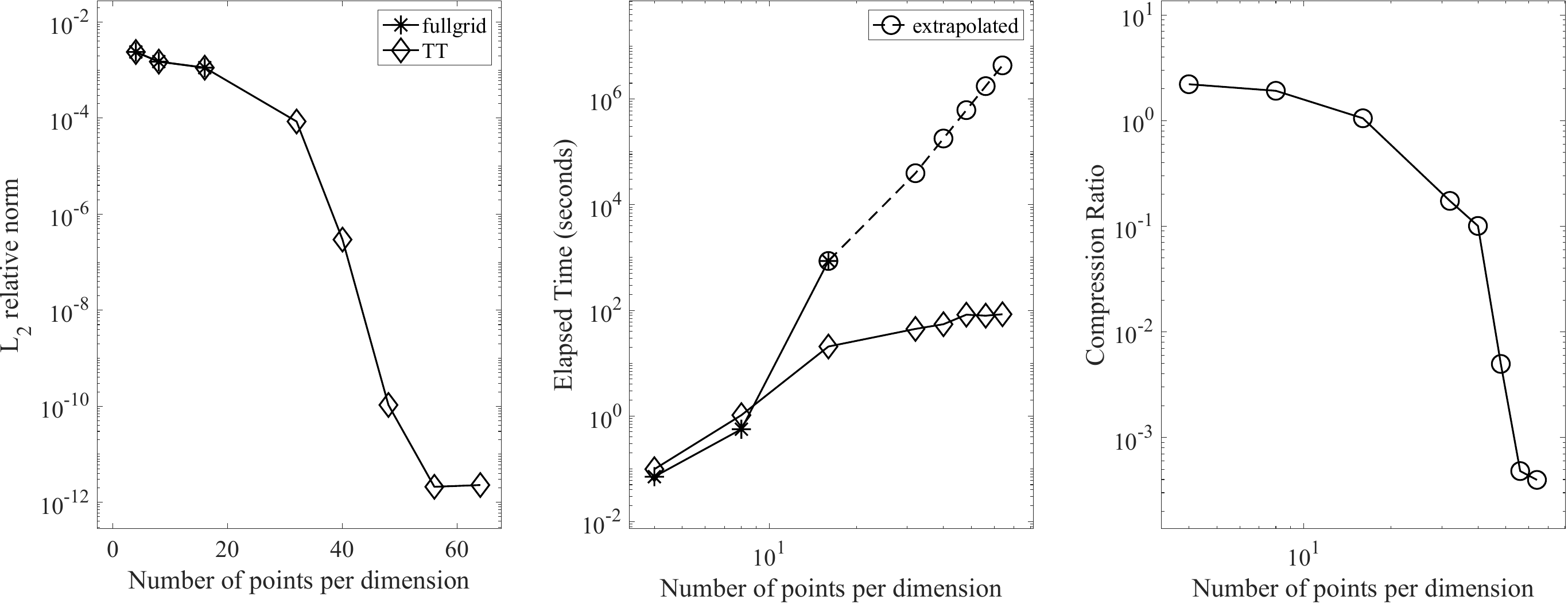}
  \caption{Test~2. Left Panel: Relative error curve in $L^2$ norm
    showing the exponential convergence of \SPSP{} schemes. Middle
    Panel: Elapsed time in seconds. \TTf{} format allows solutions in
    much higher resolution compared to the full grid format. Right
    Panel: Compression ratio of the solution. All plots are versus the
    number of points per dimension.}
  \label{VariableCDR_SPTT_SPFG}
\end{figure}

In Figure~\ref{VariableCDR_SPTT_SPFG}, we compare the \TTf{} format
and the full grid approaches, where we use \SPSP{} in both cases.
The \TTf{} approximation achieves the same level of accuracy as the
full grid computation.
The left panel plot illustrates the exponential convergence observed
in both approaches.
The error at $64$ nodes does not improve further because it reaches
the truncation tolerance of the \TTf{} format.

In the middle panel, we compare the elapsed time (in seconds) required
for both approaches.
The time taken by the full grid approach increases exponentially,
while it increases linearly in the \TTf{} format.
We conclude that the \TTf{} format provides approximately a $40$-fold
speed-up compared to the full grid computation at a grid size of $16$
nodes per dimension, and an estimated $5\times10^4$ times faster at
$64$ nodes based on extrapolated full grid time.
It is worth noting that the full grid computation is unable to handle
the grid data with 32 nodes per dimension due to memory limitations.

In the third panel, we show the compression ratio of the \TTf{} format
against the number of nodes per dimension.
As expected, the compression ratio decreases exponentially.
It is significant to highlight that the \SPSP{} scheme with \TTf{}
format achieves storage savings of approximately $10^3$ orders of
magnitude at $64$ nodes per dimension.
Additionally, Table~\ref{tab:test_2_compress} showcases the aggregated
operator $\ten{A}$ compression, demonstrating that the \TTf{} format
allows full access to terabyte-sized operators using only kilobytes of
storage.
For example, at $64$ nodes, the \TTf{} format only requires 1.81 MB of
storage to store a full grid operator with a size of 1640 terabytes.

\begin{table}[h]
  \centering
  \begin{tabular}{ |c|c|c|c|} 
    \hline
    Number of points per dimension & size of $\ten{A}$ in TB & size of $\ten{A}^{TT}$ in KB & Compress ratio  \\ 
    \hline
    8&1.66e-5  &1.74e1 & 9.73e-4\\ 
    16&1.23e-2  &9.30e1 & 7.03e-6\\ 
    32&5.10e0  &4.24e2 & 7.75e-8\\ 
    64&1.64e3  &1.81e3 & 1.03e-9\\ 
    \hline
  \end{tabular}
  \caption{Storage Cost Comparison of $\ten{A}$ between full size and
    TT-format. At 64 points, the TT-format allows compression from
    1640 TBs to about 1.81 MBs, which is about nine orders of
    magnitude.}
  \label{tab:test_2_compress}
\end{table}
\begin{figure}
  \centering
  \includegraphics[width=0.9\textwidth]{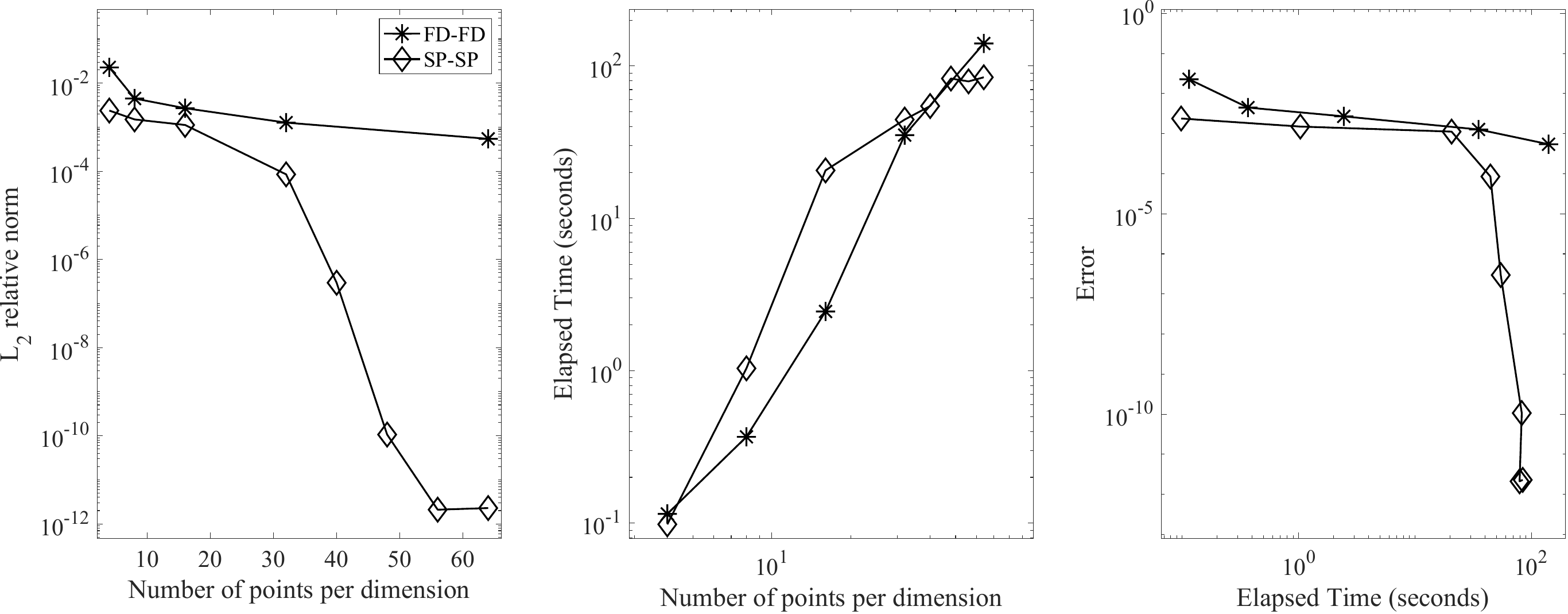}
  \caption{Test 2: Comparison between \TTf{} \SPSP{} and \TTf{}
    \FDFD{}. Left Panel: Relative error versus number of points per
    dimension. Middle Panel: Elapsed time versus number of points per
    dimension. Right Panel: Relative error versus elapsed times, which
    shows that the \TTf{} \SPSP{} is more efficient compared to \TTf{}
    \FDFD{}.}
  \label{VariableCDR_SPTT_FDTT}
\end{figure}

In Figure~\ref{VariableCDR_SPTT_FDTT}, we present a comparison between
the \FDFD{} scheme in \TTf{} format and the \SPSP{} scheme in \TTf{}
format. The left panel displays the plot of $L^2$ relative errors against the
number of nodes per dimension.
The \SPSP{} scheme with \TTf{} format achieves an approximation of the
numerical solution that is approximately accurate to the order of
$10^{-12}$, while the \FDFD{} scheme is only accurate to the order of
$10^{-3}$ at 64 nodes per dimension.

In the middle panel, we plot the elapsed times for both methods
against the number of nodes per dimension.
Both methods consume almost the same amount of time for all grid
data.

In the right panel, we plot the errors computed by both methods
against the elapsed time.
The \SPSP{} scheme generates a highly accurate numerical solution
compared to the \FDFD{} scheme.
Specifically, the \SPSP{} scheme achieves an accuracy of approximately
$10^{-12}$, while the \FDFD{} scheme only approximates the solution to
the order of $10^{-3}$ within the same time frame.
Furthermore, the accuracy increases exponentially in the case of the
\SPSP{} method.

\subsection{Test~3: Non-smooth solution}
The primary objective of this example is to validate the suitability
of the tensor-train space-time spectral collocation method for solving
the time-dependent CDR equation with less regular solutions.
To investigate this aspect, we examine the behavior of the methods on
a benchmark problem with the constant coefficients
$\kappa(t,\bold{x})=1$, $\bold{b}(t,\bold{x})=(1,1,1)$, and
$c(t,\bold{x})=1$, and the non-smooth exact solution
$u(t,x,y,z):=\sin(\pi x)\sin(\pi y)\sin(\pi z) + x^2\vert x \vert$.
As in the previous cases, computational domain is $[0,1] \times
[-1,1]^3$.

The exact solution $u$ satisfies
$\frac{d^2u(t,\bold{x})}{d\bold{x}^2}\in H^{\alpha}(0,T;L^2(\Omega))$
and $\frac{d^2u(t,\bold{x})}{d\bold{x} ^2}\in C^0(\Omega)$ for all
$t\in(0,T)$ and $\alpha>0$.
Accordingly, Theorem~\ref{Order:convergence} implies that the discrete
approximation converges quadratically.

\begin{figure}
  \centering
  \includegraphics[width=0.9\textwidth]{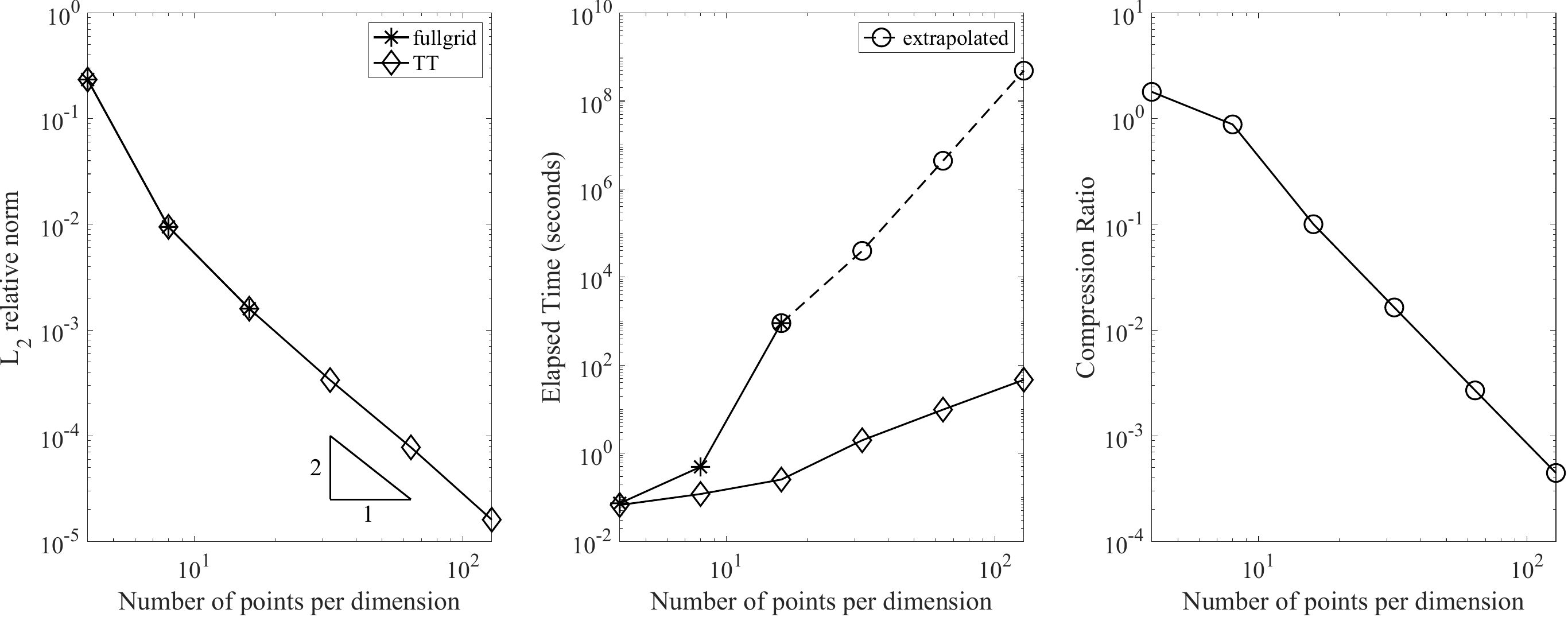}
  \caption{Test 3. Left Panel: Relative
    error curve in $L^2$ norm showing the linear convergence of
    \SPSP{} schemes, while the theoretical expected quadratic convergence is shown by the triangle. Middle Panel: Elapsed time in seconds. \TTf{}
    format allows solutions in much higher resolution compared to the
    full grid format. Right Panel: Compression ratio of the
    solution. All plots are versus the number of points per
    dimension.}
  \label{fig:test3-FG-TT}
\end{figure}

Indeed, we observe such a convergence order is clear in the left panel
of Fig.~\ref{fig:test3-FG-TT}.
Due to the lower regularity of the exact solution, a highly refined
mesh is required to achieve high accuracy, which is only feasible with
the \TTf{} format.

In the middle panel, the plot shows that the \SPSP{} method is
approximately $3500$ times faster than the full-grid \SPSP{} method
using $16$ nodes, and about $1\times10^7$ times faster using $128$
nodes, this latter being extrapolated on the full-grid.
The compression ratio is approximately $5\times10^{-4}$.
Since the coefficient functions are the same as in Test Case 1, the
compression of the operators in this test case is similar to the ones
shown in Table \ref{tab:test_1_compress}.

\begin{figure}
  \centering
  \includegraphics[width=0.9\textwidth]{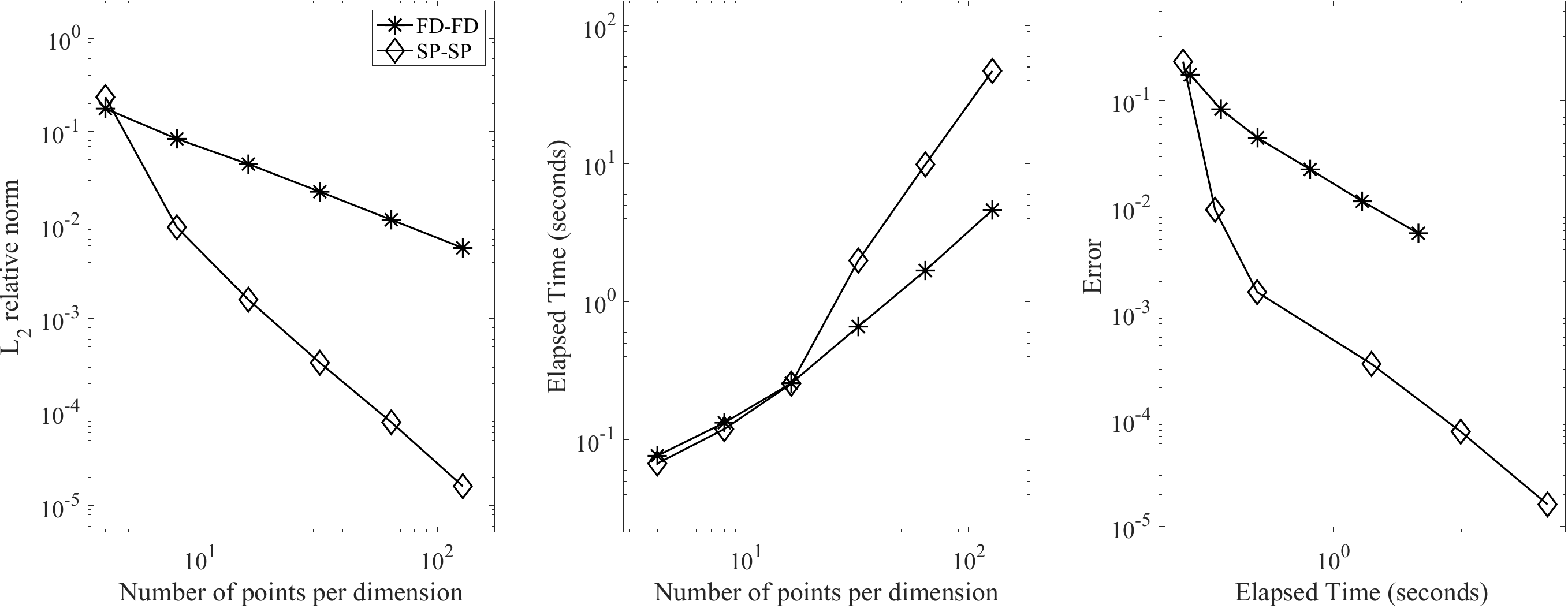}
  \caption{Test case-3:Comparison between \TTf{} \SPSP{} and \TTf{}
    \FDFD{}. Left Panel: Relative error versus number of points per
    dimension. Middle Panel: Elapsed time versus number of points per
    dimension. As the grid becomes larger, the TT \SPSP{} requires more computational time.
    Right Panel: Relative error versus elapsed times, which
    shows that the \TTf{} \SPSP{} is more efficient compared to \TTf{}
    \FDFD{}.}
  \label{fig:test3-TT}
\end{figure}

Fig.~\ref{fig:test3-TT} compares the \TTf{} \FDFD{} and \TTf{} \SPSP{}
methods.
The \TTf{} \SPSP{} method exhibits linear convergence compared to
\TTf{} \FDFD{} with smaller approximation errors.
Furthermore, as the grid size increases, the TT \SPSP{} method
requires more elapsed time.
The right panel, which plots the error versus elapsed time, clearly
shows that the \TTf{} \SPSP{} method remains more efficient than the
\TTf{} \FDFD{} method for this test case with a less regular solution.

\section{Conclusions}
\label{conclusions}
In this work, we introduce a method called the TT spectral collocation space-time (TT-SCST) approach for solving time-dependent convection-diffusion-reaction equations. Time is considered as an additional dimension, and the spectral collocation technique is applied in both space and time dimensions. This space-time approach demonstrates exponential convergence property for smooth functions. However, it involves solving a large linear system for solutions at all time points. Utilizing the tensorization process, we convert the linear system into TT-format and subsequently solve the TT linear system. Our numerical experiments validate that the TT-SCST method achieves significant compression of terabyte-sized matrices to kilobytes, leading to a computational speedup of tens of thousands of times, while maintaining the same level of accuracy as the full-grid space-time scheme.


\section*{Acknowledgments}
The authors gratefully acknowledge the support of the Laboratory
Directed Research and Development (LDRD) program of Los Alamos
National Laboratory under project number 20230067DR.
Los Alamos National Laboratory is operated by Triad National Security,
LLC, for the National Nuclear Security Administration of
U.S. Department of Energy (Contract No.\ 89233218CNA000001).


\bibliographystyle{plain}
\bibliography{ttrain}

\setcounter{section}{0}
\appendix
\section{Notation, basic definitions, and operations with tensors}
\label{Appendix_B}
Let $d$ be a positive integer.
A $d$-dimensional tensor $\ten{A}\in\mathbb{R}^{n_1\times
  n_2\times\ldots\times n_d}$ is a multi-dimensional array with $d$
indices and $n_k$ elements in the $k$-th direction for
$k=1,2,\ldots,d$.
We say that the number of dimensions $d$ is the \emph{order of the
tensor}.
\subsection{Kronecker product}
The Kronecker product $\bigotimes$ of matrix
$\mat{A}=(a_{ij})\in\mathbb{R}^{m_A \times n_A}$ and matrix
$\mat{B}=(b_{ij})\in\mathbb{R}^{m_B \times n_B}$ is the matrix
$\mat{A}\otimes \mat{B}$ of size
$N_{\mat{A}\otimes\mat{B}}=(m_Am_B)\times(n_An_B)$ defined as:
\begin{equation}
  \mat{A}\otimes\mat{B} = 
  \begin{bmatrix}
    a_{11}\mat{B} &a_{12}\mat{B} &\cdots & a_{1n_A}\mat{B}\\
    a_{21}\mat{B} &a_{22}\mat{B} &\cdots & a_{2n_A}\mat{B}\\
    \vdots & \vdots & \ddots & \vdots \\
    a_{m_A1}\mat{B} &a_{m_A2}\mat{B} &\cdots & a_{m_An_A}\mat{B}\\
  \end{bmatrix}.
  \label{def:kronecker}
\end{equation}
Equivalently, it holds that
$\big(\mat{A}\otimes\mat{B}\big)_{ij}=a_{i_Aj_A}b_{i_Bj_B}$, where
$i=i_B+(i_A-1)m_B$, $j=j_B+(j_A-1)m_B$, with $i_A=1,\ldots,m_A$,
$j_A=1,\ldots,n_A$, $i_B=1,\ldots,m_B$, and $j_B=1,\ldots,n_B$.

\subsection{Tensor product}
There is a relation between Kronecker product and tensor product: the
Kronecker product is a \textit{special} bilinear map on a pair of
vector spaces consisting of matrices of a given dimensions (it
requires a choice of basis), while the tensor product is a
\textit{universal} bilinear map on a pair of vector spaces of any
sort, hence, it is more general.

Here we define the tensor product of two vectors
$\mat{a}=(a_i)\in\mathbb{R}^{n_A}$ and
$\mat{b}=(b_i)\in\mathbb{R}^{n_B}$, which produces the
\textbf{two dimensional-tensor}
$\big(\mat{a}\tenprod\mat{b}\big)$ of size
$N_{\mat{a}\tenprod\mat{b}}=n_A\times n_B$ defined as:
\begin{equation}
  \big(\mat{a}\tenprod\mat{b}\big)_{ij} = a_{i}b_{j}
  \qquad
  i=1,2,\ldots,n_A,\,
  j=1,2,\ldots,n_B.
  \label{def:outer:product:vectors}
\end{equation}
Note that
$\mat{a}\tenprod\mat{b}=\mat{a}\otimes\mat{b}^T$.
Similarly, the tensor product of matrix
$\mat{A}=(a_{ij})\in\mathbb{R}^{m_A \times n_A}$ and matrix
$\mat{B}=(b_{kl})\in\mathbb{R}^{m_B \times n_B}$ produces the
four-dimensional tensor of size $N_{\mat{A}\tenprod\mat{B}}=m_A \times
n_A \times m_B \times n_B$, with elements:
\begin{equation}
  \big(\mat{A}\tenprod\mat{B}\big)_{ijkl} = a_{ij}b_{kl},
  \label{def:outer_product:matrices}
\end{equation}
for
$i=1,2,\ldots,m_A$,
$j=1,2,\ldots,n_A$,
$k=1,2,\ldots,m_B$,
$l=1,2,\ldots,n_B$.

\section{Construction of \texorpdfstring{$\ten{A}^{map,TT} \text{ and } 
\label{APP:A_map_tt}
\ten{G}^{TT}$}{}} Here we provide details about the construction of
$\ten{A}^{map,TT}$ and $\ten{G}^{Bd,TT}$.

\begin{equation*}
  \ten{A}^{map,TT} = \sum_{k = t,D,C,R} \ten{A}_k^{map,TT},
\end{equation*}
where 
\begin{equation}
\begin{split}
  &\ten{A}^{map,TT}_t = \bold{S}_t(\ten{I}_t,\ten{I}_t)\tenprod \bold{I}^{N+1}(\ten{I}_s,:)\tenprod \bold{I}^{N+1}(\ten{I}_s,:)\tenprod \bold{I}^{N+1}(\ten{I}_s,:),\\
\end{split}
\end{equation}
\begin{equation}
  \begin{split}
    &\Delta^{map,TT}=
    \begin{aligned}[t]
      &\bold {I}^{N+1}(\ten{I}_t,:)\tenprod \bold{S}_{xx}(\ten{I}_s,:) \tenprod \bold {I}^{N+1}(\ten{I}_s,:)\tenprod \bold {I}^{N+1}(\ten{I}_s,:)\\
      +&\bold {I}^{N+1}(\ten{I}_t,:) \tenprod \bold{I}^{N+1}(\ten{I}_s,:) \tenprod \bold{S}_{yy}(\ten{I}_s,:)  \tenprod\bold{I}^{N+1}(\ten{I}_s,:)\\
      +&\bold {I}^{N+1}(\ten{I}_t,:) \tenprod \bold{I}^{N+1(\ten{I}_s,:} \tenprod  \bold{I}^{N+1(\ten{I}_s,:} \tenprod\bold{S}_{zz}(\ten{I}_s,:),
    \end{aligned}\\
    &\ten{A}^{map,TT}_D = \ten{K}^{TT,op}\Delta^{map,TT}
  \end{split}
\end{equation}

\begin{equation}
  \begin{split}
    &\nabla_x^{map,TT} = \bold {I}^{N+1}(\ten{I}_t,:)\tenprod \bold{S}_{x}(\ten{I}_s,:) \tenprod \bold{I}^{N+1}(\ten{I}_s,:)\tenprod \bold{I}^{N+1}(\ten{I}_s,:) \\
    &\nabla_y^{map,TT} = \bold {I}^{N+1}(\ten{I}_t,:)\tenprod \bold{I}^{N+1}(\ten{I}_s,:) \tenprod\bold{S}_{y}(\ten{I}_s,:) \tenprod \bold{I}^{N+1}(\ten{I}_s,:) \\
    &\nabla_z^{map,TT} = \bold {I}^{N+1}(\ten{I}_t,:)\tenprod \bold{I}^{N+1}(\ten{I}_s,:) \tenprod \bold{I}^{N+1}(\ten{I}_s,:) \tenprod\bold{S}_{z}(\ten{I}_s,:)  \\
    &\ten{A}_C^{map,TT} = \ten{B}^{x,TT,op}\nabla_x^{map,TT} + \ten{B}^{y,TT,op}\nabla_y^{map,TT} + \ten{B}^{z,TT,op}\nabla_z^{map,TT}
  \end{split}
\end{equation}

\begin{equation}
    \ten{A}_R^{map,TT} = \ten{C}^{TT,op}(\bold {I}^{N+1}(\ten{I}_t,:)\tenprod \bold{I}^{N+1}(\ten{I}_s,:) \tenprod \bold{I}^{N+1}(\ten{I}_s,:))
\end{equation}

Next, we show how to construct the $\ten{G}^{Bd}$ tensor.
The tensor $\ten{G}^{Bd}$ is a $(N+1)\times (N+1) \times (N+1) \times
(N+1)$, in which only the BC/IC elements are computed. Other elements
are zeros. The TT tensor $\ten{G}^{Bd}$ is constructed using the cross
interpolation.

Then, the boundary tensor $\ten{F}^{Bd,TT}$ is computed as:
\begin{equation}
  \ten{F}^{Bd,TT}= \ten{A}^{map,TT}\ten{G}^{Bd,TT}
\end{equation}

\end{document}
